\documentclass[11pt]{amsart}
\usepackage[T1]{fontenc}
\usepackage{lmodern}
\usepackage{needspace}
\usepackage[hidelinks]{hyperref}
\newtheorem{theorem}{Theorem}[section]
\newtheorem{lemma}[theorem]{Lemma}
\newtheorem{corollary}[theorem]{Corollary}
\newcommand{\Inn}{\operatorname{Inn}}
\newcommand{\Sym}{\operatorname{Sym}}

\title[Almost all finite racks and quandles are 2-reductive]{Almost all finite racks and quandles are 2-reductive}
\author{Amir Saki}
\address{Universit\'e du Qu\'ebec \`a Trois-Rivi\`eres, Trois-Rivi\`eres, Qu\'ebec, Canada}
\email{amir.saki.math@gmail.com}
\author{Usef Faghihi}
\address{Universit\'e du Qu\'ebec \`a Trois-Rivi\`eres, Trois-Rivi\`eres, Qu\'ebec, Canada}
\email{usef.faghihi@uqtr.ca}

\begin{document}
\begin{abstract}
Vojt\v{e}chovsk\'y and Yang asked whether the proportions of \(2\)-reductive racks and quandles of order \(n\), counted up to isomorphism, have limits as \(n\) tends to infinity. We prove that both limits are \(1\). More precisely, for some constant \(c>0\) and all sufficiently large \(n\), the proportion of racks of order \(n\) whose inner automorphism group has an orbit of size at least three is at most \(2^{-cn}\). The same bound holds for quandles.
\end{abstract}
\maketitle

\section{Introduction}

Let \(R\) be a set with a binary operation \(*\). For each \(x\in R\), let \(L_x\) denote the map \(y\mapsto x*y\). We say that \(R\) is a rack if each \(L_x\) is a bijection and
\[
x*(y*z)=(x*y)*(x*z)
\]
for all \(x,y,z\in R\). A rack is a quandle if \(x*x=x\) for each \(x\in R\). We write \(\Inn(R)=\langle L_x:x\in R\rangle\). For each \(g\in\Inn(R)\) and \(x\in R\), we have
\begin{equation}\label{eq:equivariance}
L_{g(x)}=gL_xg^{-1}.
\end{equation}
A rack is \(2\)-reductive if \((x*y)*z=y*z\) for all \(x,y,z\in R\). By \eqref{eq:equivariance}, this is equivalent to \(\Inn(R)\) being abelian.

Let \(r(n)\) and \(q(n)\) denote the numbers of isomorphism classes of racks and quandles of order \(n\), respectively. Let \(r_{\mathrm{2\text{-}red}}(n)\) and \(q_{\mathrm{2\text{-}red}}(n)\) denote the corresponding numbers of \(2\)-reductive racks and quandles. Vojt\v{e}chovsk\'y and Yang asked whether the ratios \(r_{\mathrm{2\text{-}red}}(n)/r(n)\) and \(q_{\mathrm{2\text{-}red}}(n)/q(n)\) have limits, and what those limits are \cite[Problem~6.2]{VY}.

Blackburn constructed a large family of racks and quandles whose inner automorphism groups have orbits of size at most two \cite{Blackburn}. Ashford and Riordan determined the leading term of \(\log_2 r(n)\) \cite[Theorem~1.1]{AR}. Ashford proved that if a sequence of prescribed orbit partitions supports \(2^{(1/4+o(1))n^2}\) labelled racks on \(\{1,\ldots,n\}\), then \(n-o(n)\) points lie in two-element blocks \cite[Proposition~7.4.4]{Ashford}. To determine the ratios above, we count racks with a larger orbit even when only a few elements lie outside two-element orbits.

\Needspace{5\baselineskip}
\begin{theorem}\label{thm:main}
There are constants \(c>0\) and \(n_0\) with the following property. For each \(n\ge n_0\), the proportion of isomorphism classes of racks of order \(n\) whose inner automorphism group has an orbit of size at least three is at most \(2^{-cn}\). The same conclusion holds for quandles.
\end{theorem}

\begin{corollary}\label{cor:problem}
Both limits in Vojt\v{e}chovsk\'y and Yang's Problem~6.2 are \(1\). Moreover, the proportion of isomorphism classes whose inner automorphism group is an elementary abelian \(2\)-group tends to \(1\), for both racks and quandles.
\end{corollary}

\begin{proof}
If all orbits of \(\Inn(R)\) have size at most two, then restriction to the two-element orbits embeds \(\Inn(R)\) in a product of copies of \(C_2\). Thus \(\Inn(R)\) is an elementary abelian \(2\)-group, and \(R\) is \(2\)-reductive. Apply Theorem~\ref{thm:main}.
\end{proof}

\section{Restrictions between two orbits}

\begin{lemma}\label{lem:pair}
Let \(R\) be a finite rack, and let \(P=\{u,v\}\) and \(O\) be orbits of \(G=\Inn(R)\), where \(|O|=d\ge3\). Assume that the permutation group \(H=G|_O\) and all restrictions \(L_x|_O\), for \(x\in O\), are given. There is an absolute constant \(\delta>0\) such that the number of possibilities for
\[
\bigl(L_u|_O,L_v|_O,(L_x|_P)_{x\in O}\bigr)
\]
is at most \(2^{(1-\delta)d}\).
\end{lemma}

\begin{proof}
Put \(a=L_u|_O\) and \(b=L_v|_O\). By \eqref{eq:equivariance}, the conjugacy class of \(a\) in \(H\) is \(\{a,b\}\). Thus \(a\) determines \(b\). The derived subgroup \(H'\) acts trivially on this class, so \(a\in C_H(H')\).

By \eqref{eq:equivariance} and the commutativity of \(\Sym(P)\), the restriction \(L_x|_P\) is constant as \(x\) ranges over \(O\). If \(a\ne b\), then \(L_x|_P\) interchanges \(u\) and \(v\) exactly when \(L_x|_O\) conjugates \(a\) to \(b\). Hence there is at most one choice of the common restriction if \(a\ne b\), and at most two choices if \(a=b\). The latter case occurs only for \(a\in Z(H)\). Therefore the number of possibilities is at most
\[
|C_H(H')|+|Z(H)|.
\]

Assume that \(H'\) has \(t\) orbits of size \(s\) on \(O\). These sizes are equal because \(H'\) is normal in the transitive group \(H\), and \(st=d\). The image of \(H\) on the set of \(H'\)-orbits is abelian and transitive. It therefore has order \(t\). Given the permutation induced by an element of \(C_H(H')\) on these orbits, the element is determined by the image of one point from each orbit. Accordingly, \(|C_H(H')|\le ts^t\). The center \(Z(H)\) acts freely on \(O\), since an element of the center fixing one point fixes every point. Thus \(|Z(H)|\le d\).

For positive integers \(s,t\), the inequalities \(t\le 2^{t-1}\) and \(s\le 2^{s-1}\) give \(ts^t\le 2^{d-1}\). Since \(d<2^{d-1}\) for \(d\ge3\), we obtain \(ts^t+d<2^d\). Also,
\[
\log_2(ts^t)\le \log_2 d+\frac{\log_2 3}{3}d,
\]
since \((\log_2 s)/s\le(\log_2 3)/3\) for all positive integers \(s\). These bounds give one constant \(\delta>0\) that works for all \(d\ge3\).
\end{proof}

\section{Racks with few elements outside two-element orbits}

Suppose that \(R\) has \(m\) two-element orbits and \(\ell\) singleton orbits. Let \(B\) be the union of its orbits of size at least three and put \(k=|B|\). Then \(n=2m+\ell+k\). We consider \(k\ge3\) and \(m\ge2\).

The translations on \(R\setminus B\), restricted to \(R\setminus B\), are described by a binary \(m\times m\) matrix and a multiset of \(\ell\) binary vectors of length \(m\). The matrix entry in position \((i,j)\) records whether a translation from pair \(i\) interchanges the elements of pair \(j\). This entry is well defined by \eqref{eq:equivariance}, since \(\Sym(2)\) is abelian. Each singleton gives one vector recording its action on the pairs. We allow zero columns because a pair might be interchanged only by translations from \(B\). For quandles the diagonal entries are zero. Two such data sets are equivalent if a simultaneous permutation of the matrix rows and columns, applied also to the vector coordinates, takes one to the other. Let \(D_{m,\ell}\) denote the number of equivalence classes, using the applicable diagonal convention.

\begin{lemma}\label{lem:rigid}
As \(m\) tends to infinity, the proportion of these equivalence classes whose matrices have no zero column and are fixed by no nonidentity simultaneous row and column permutation tends to \(1\), uniformly in \(\ell\).
\end{lemma}

\begin{proof}
On indexed matrices the fraction with a zero column is at most \(m2^{-(m-1)}\). Let \(\sigma\in\Sym(m)\) move \(r\) indices. It moves at least \(r(m-1)\) available matrix positions, since every position in a moved row is moved. Each cycle of moved positions has at least two positions. Thus the fraction of matrices fixed by \(\sigma\) is at most \(2^{-r(m-1)/2}\). The sum of these fractions over nonidentity permutations is at most
\[
\sum_{r=2}^m m^r 2^{-r(m-1)/2}=o(1).
\]
In particular, a \(1-o(1)\) fraction of indexed matrices have both required properties.

Let \(N\) be the number of indexed matrices and \(M\) the number of multisets of \(\ell\) vectors on the indexed coordinates. A permutation fixes at most \(M\) times as many decorated matrices as it fixes matrices. Burnside's lemma therefore gives \(D_{m,\ell}\le (1+o(1))NM/m!\). The classes with both required properties have \(m!\) representatives each, so their number is at least \((1-o(1))NM/m!\). The error bounds are independent of \(M\).
\end{proof}

Write \(T(n)=r(n)\) in the rack case and \(T(n)=q(n)\) in the quandle case. Let \(E(m,\ell,k)\) denote the number of isomorphism classes with these orbit counts.

\begin{lemma}\label{lem:ratio}
There is an absolute constant \(C\) such that, for all sufficiently large \(m\),
\[
\frac{E(m,\ell,k)}{T(n)}
\le 2^{-\delta mk+C(k^2+k\ell+k\log_2 n)}.
\]
\end{lemma}

\begin{proof}
Fix a representative of one of the \(D_{m,\ell}\) classes. On a labelled copy of \(B\), choose the restricted rack operation and the permutation group \(G|_B\), where \(G=\Inn(R)\). There are \(2^{O(k^2)}\) choices in total. The bound for the operation follows from Ashford and Riordan \cite[Theorem~1.1]{AR}, in its labelled form. The group bound follows from Pyber's estimate for the number of subgroups of \(\Sym(k)\) \cite{Pyber}, as stated in \cite[Theorem~4.2.3]{Ashford}.

The group \(G|_B\) determines the orbits in \(B\). For each singleton \(\{y\}\) of \(R\) and each orbit \(O\) in \(B\), \eqref{eq:equivariance} gives \(L_y|_O\in Z(G|_O)\), so there are at most \(|O|\) choices. This gives at most \(2^{O(k\ell)}\) choices for the restrictions from singletons to \(B\). Every translation fixes the singleton points. By Lemma~\ref{lem:pair}, the remaining restrictions between pairs and \(B\) give at most \(2^{(1-\delta)mk}\) choices. Hence
\begin{equation}\label{eq:upper}
E(m,\ell,k)\le D_{m,\ell}\,2^{(1-\delta)mk+O(k^2+k\ell)}.
\end{equation}

For a lower bound on \(T(n)\), take only classes from Lemma~\ref{lem:rigid} that satisfy both properties. For each such class, keep a fixed indexed representative and add \(h=\lfloor k/2\rfloor\) ordered pairs and, if \(k\) is odd, a distinguished singleton. Make each translation constant as a function of its source point within each block, and let it act by independent interchanges on target pairs. Keep the old data. Among the \(2mh\) entries between old and new pairs, set the \(h\) entries from the first old pair to the new pairs to \(1\). Choose the remaining \(2mh-h\) entries freely. If \(k\) is odd, also choose freely the \(m\) entries from the distinguished singleton to the old pairs. Set all other new entries to zero. The translations commute and preserve the blocks, so \(L_{L_x(y)}=L_xL_yL_x^{-1}\) for each \(x,y\). Hence the operations are racks. They are quandles when the old data have zero diagonal. Each pair is an orbit, by the nonzero columns for old pairs and the fixed bits for new pairs.

There are \(2^{mk-h}\) choices per class of old data. Given an output isomorphism class, there are at most \(n^{h+1}\) ways to identify its ordered new pairs and, when present, its distinguished singleton. Once these are identified, restriction recovers the old data class. Its matrix has no nonidentity index symmetry, so an isomorphism with the chosen indexed representative determines the old pair indices uniquely. Interchanges within individual pairs and permutations of old singletons leave all free bits unchanged. Thus each identification recovers at most one choice of the free bits. By Lemma~\ref{lem:rigid}, for all sufficiently large \(m\),
\[
T(n)\ge D_{m,\ell}\,2^{mk-O(k\log_2 n)}.
\]
Together with \eqref{eq:upper}, this proves the lemma.
\end{proof}

\section{Proof of Theorem~\ref{thm:main}}

Choose a fixed \(\varepsilon>0\) sufficiently small. If \(k+\ell\le\varepsilon n\) and \(k\ge3\), then \(m\ge (1-\varepsilon)n/2\). By Lemma~\ref{lem:ratio}, for a constant \(c_1>0\) and all sufficiently large \(n\), the proportion with fixed \(m,\ell,k\) is at most \(2^{-c_1nk}\). There are at most \((n+1)^2\) triples to consider. Thus this range contributes at most \(2^{-c_2 n}\) for some \(c_2>0\).

For the other range, let \(\mathcal P\) be a partition of \(\{1,\ldots,n\}\), and let \(F(\mathcal P)\) be the number of racks on this set whose orbit partition is \(\mathcal P\). Ashford's Proposition~7.4.4 implies that, for each fixed \(\varepsilon>0\), there is an \(\eta>0\) such that
\[
F(\mathcal P)\le 2^{(1/4-\eta)n^2}
\]
for all sufficiently large \(n\) whenever at least \(\varepsilon n\) points of \(\mathcal P\) lie outside its two-element parts. To see the uniform implication, suppose that no such \(\eta\) exists and choose a sequence of counterexamples along increasing orders whose normalized logarithmic counts tend to \(1/4\). The upper bound of Ashford and Riordan ensures that these counts cannot have a larger limit. Fill the other orders with partitions into pairs and at most one singleton. The counts for these partitions also have normalized logarithmic limit \(1/4\), by the binary matrix construction with the diagonal zero and all columns nonzero. Proposition~7.4.4 then says that the number of points outside pairs is \(o(n)\), a contradiction.

There are at most \(n^n\) partitions. The total number of isomorphism classes in either case satisfies \(T(n)\ge 2^{n^2/4-O(n\log n)}\), by Blackburn's construction \cite{Blackburn}. Hence the range \(k+\ell>\varepsilon n\) contributes at most \(2^{-c_3n^2}\) for some \(c_3>0\). Combining the two ranges proves Theorem~\ref{thm:main}.

\section*{Use of generative AI}
Under the authors' direction and supervision, OpenAI's ChatGPT was used in parts of the development and verification of the arguments and in the preparation of this manuscript. The authors take responsibility for the contents of the manuscript.


\begin{thebibliography}{9}
\bibitem{Ashford}
M. Ashford, \emph{Graphs of algebraic objects}, DPhil thesis, University of Oxford, 2016.
\href{https://ora.ox.ac.uk/objects/uuid:a666c588-7653-4e1c-a60d-5841f1019a1e}{Oxford University Research Archive}.

\bibitem{AR}
M. Ashford and O. Riordan, \emph{Counting racks of order \(n\)}, Electron. J. Combin. \textbf{24} (2017), no.~2, Paper 2.32, 20 pp.
\href{https://doi.org/10.37236/6330}{doi:10.37236/6330}.

\bibitem{Blackburn}
S. R. Blackburn, \emph{Enumerating finite racks, quandles and kei}, Electron. J. Combin. \textbf{20} (2013), no.~3, Paper 43, 9 pp.
\href{https://doi.org/10.37236/3262}{doi:10.37236/3262}.

\bibitem{Pyber}
L. Pyber, \emph{Enumerating finite groups of given order}, Ann. of Math. (2) \textbf{137} (1993), no.~1, 203--220.

\bibitem{VY}
P. Vojt\v{e}chovsk\'y and S. Y. Yang, \emph{Enumeration of racks and quandles up to isomorphism}, Math. Comp. \textbf{88} (2019), no.~319, 2523--2540.
\href{https://doi.org/10.1090/mcom/3409}{doi:10.1090/mcom/3409}.
\end{thebibliography}
\end{document}